%% file: orderunitcancellationinRK
\input amstex

\input generic_macros

\input papermacros_entirefunctions

\def\Z{\text{\bf Z}}


\def\ti#1{\tilde{#1}}

\let\exist=\exists

\def\Z{\text{\bf Z}}

\long\def\subsubtitle #1.{\vskip 3pt \noindent {\it #1 }}

\long\def\subtitle #1\par#2{\vskip 3pt \noindent {\it #1}\par \noindent
#2}

\def\ti#1{\tilde{#1}}

\def\ideal(#1){\(#1\)}
\def\oid(#1){< #1 >}
\def\id#1{< #1 >}

\def\FF{{\Cal F}}
\def\RR{{\Cal R}}
\def\T{\pmb T}

\Title In praise of order units

\Abstract We show that the ordered rings naturally associated to compact
convex polyhedra with interior satisfy a positivity property known as
order unit cancellation, and obtain other general positivity results as
well. 

\vskip4pt 
\noindent {\it David Handelman}\plainfootnote{$^1$}{Supported in part
by a Discovery Grant from NSERC.}\vskip 4pt 

\noindent An element $u$ of a partially ordered abelian group $G$ is
called an {\it order unit\/} (relative to $G$) if for all $g$ in $G$,
there exists a positive integer $M$ \st $-Mu \leq g \leq Mu$. A {\it
trace\/} of $G$ is a positive additive nonzero function $G \to \R$,
and a {\it pure trace\/} is a trace that is not a nontrivial convex
linear combination of other traces. The partially ordered abelian group
is {\it unperforated\/} if for  $g$ in $G$ and positive integer $n$, $ng
\in G^+$ implies $g \in G^+$. When
$G$ is unperforated and admits an order unit, then $u$ will be an order
unit if and only if for all (pure) traces $\tau$ on $G$, $\tau (g) > 0$
[GH, Lemma 4.1]. An {\it order ideal\/} $I$ of a partially ordered abelian
group $G$ is a subgroup with two additional properties, {\it convexity\/}
($0 \leq g \leq h$, $h
\in I$, $g \in G$ entails $g \in I$) and directedness ($I$ is spanned by
its positive elements).

In this paper, a {\it partially ordered ring\/} $R$ will be a commutative
unital ring with a subset $R^+$ making $R$ into a partially ordered group
with respect to addition, \st in addition, $R^+ \cdot R^+ = R^+$, $R$ is
unperforated, and $1$ is an order unit (relative to $R$). By [H85, 1.2(a)], any
order ideal of a partially ordered ring is itself an ideal (this depends
of course on $1$ being an order unit). The pure traces on an ordered
ring with $1$ as order unit, normalized so that $1 \to 1$, are exactly
the multiplicative traces. The set of these, equipped with the
point-open topology, form a compact set (the extremal boundary of the
Bauer simplex consisting of all the normalized traces).  

We say the partially ordered ring $R$ satisfies {\it order unit
cancellation\/} (informally described as {\it order cancellation\/} in [H87A, p\,22]) if whenever $u$ is an order unit of $R$ and a nonzero
divisor, and $a$ is an element of $R$, $ua \geq 0$ implies $a \geq 0$.
(As usual, we use $r \geq s$ to denote $r-s \in R^+$.) This is a
relatively weak property for an ordered ring, but it has been
useful. For example, it is the key tool in  [H87A, II.1], showing that if
$P$, $Q$, and $f$ are real polynomials in several variables \st $P$ and
$Q$ have no nonnegative coefficients and the same monomials appear, then 
$$\eqalign{
\text{there exists $m$ \st $P^m f$ }& \text{has no negative coefficients
} 
\cr
 \text{if and only if }&\text{there exists $M$ \st $Q^M f$ has no negative
coefficients.}\cr }$$ (The actual result is sharper  than this.)

If $U \equiv U(R)$ denotes the set of nonzero divisor order units of the
partially ordered ring $R$, then $U$ is multiplicatively closed (and also
additively closed if $R$ is a domain, but this seems to play a limited
role in what follows). We may thus form the overring $U^{-1}R$, which can
be equipped with the direct limit partial ordering, $(U^{-1}R )^+$,
obtained by taking the directed limit over all maps of the form $\Arrow
\times u; R.R$, every element $u$ appearing infinitely often. Then 
$$ (U^{-1}R )^+ = \Set{u^{-1}r}{u \in U, \text{there exists $v$ in $U$
\st $vr \in R^+$}}. 
$$ 
(We could just as well take this as the definition of $(U^{-1}R )^+$
rather than using the limit formulation; however, direct limits appear
throughout the examples.)  Then order unit cancellation is equivalent to
$(U^{-1}R )^+ \cap R = R^+$. The effect of this localization is to
eliminate  the multiplicative homomorphisms that are not
traces (i.e., not positive), so that the maximal ideal space of $U^{-1}R$ is identifiable in the obvious way with the pure trace space; this is used in [H87A, App A] to deduce
properties of decompositions of convex polyhedra as sums of others. This
should also be the precursor to the use of localization at prime order
ideals. 

Our principle result will be to show that a class of partially ordered
rings studied in [H88]  satisfies order unit
cancellation. Let $K$ be a compact convex polyhedron (a \quotes{\it
polytope\/}) with interior, in $\R^n$. For each facet, $F_i$ of $K$,
choose a linear form (unique up to positive scalar multiple) $\beta_i =
\sum a_i x_i + a_0$ (where $a_0, a_1, \dots, a_n$ are real) \st
$\beta_i|F_i \equiv 0$ and $\beta_i |K \geq 0$. Let $R = \R[K]$ denote the
polynomial ring $\R[x_1, \dots, x_n]$ equipped with the positive cone
generated additively and multiplicatively (allowing positive scalars) by
$\brcs{\beta_i}$. It was shown in [H88, I.2] that this is a partially ordered
ring in the sense used here, including the fact that $1$ is an order
unit. We will show that $\R[K]$ has order unit cancellation. We will also
obtain a number of positivity results in a similar vein to order unit
cancellation. 

Before doing so, we note a number of classes of examples and
non-examples. 

\subsubtitle Trivial examples. Let $X$ be a compact set, and let $R$ be
any unital subring of the ring of continuous real-valued functions,
$C(X)$, equipped with the pointwise order (that is, for $f,g$ in $R$,
$f-g \in R^+$ if and only if $f(x) \geq g(x)$ for all $x $ in $X$). Order
unit cancellation is obvious.

A partially ordered abelian group is {\it simple\/} if it has no proper
order ideals, or what amounts to the same thing, all nonzero positive
elements are order units. If the partially ordered ring $R$ is simple as
a partially ordered abelian group, then obviously order unit cancellation
holds. 

If $F$ is a formally real subfield of the algebraic closure of the
rationals equipped with the sum of squares ordering, then $F$ is a
partially ordered ring that is simple as a partially ordered group---but
this falls in the previous class of examples, with $X$ being the Galois
group of the real algebraic closure of the rationals. 

A larger class consists of partially ordered rings $R$ (with $1$)
equipped with a unital algebra map $\Arrow \pi;R . C(X)$ (not necessarily
one to one), \st a nonzero element of $R$ is in $R^+$ if and only its
image in $C(X)$ is strictly positive (as a function). This in fact
describes completely the class of partially ordered rings which are
simple as partially ordered groups. 
\qed

Before dealing with more examples, we cite a result that we will prove
later. 

\Lem Theorem 1. Let $R$ be a noetherian partially ordered ring
(unperforated, and with $1$ as order unit). Suppose that for every proper
order ideal
$I$ of
$R$ the following holds:
\itemitem{} for every prime ideal $P$ minimal  over $I$, there exists a
proper order ideal $J$ of $R$ and an integer $k$ \st $P^k \subseteq
J$.{\par} \noindent Then $R$ satisfies order unit cancellation.

\subsubtitle Toy examples. Let $R$ be the ring $\R[x]$, which we equip
with two partial orderings, but these can obviously be modified. 
$$\eqalign{ 
R_1^+ & = \Set{x^j(1+x) f}{j \geq 1, \ f|[0,1] > 0}  \cup \Set{f}{f|[0,1] > 0}\cr 
R_2^+ & = \Set{x^j f}{j = 0,2, 3, \dots; \ f|[0,1] > 0} 
}$$ In both cases, we
see the putative positive cones are closed under addition, positive
scalar multiplication, and multiplication, and $1$ is an order unit. In
both cases, the pure traces are point evaluations at the points of
$[0,1]$, so the order units are precisely the polynomials that are
strictly positive thereon. 

In the case of $R_1$, we note that $x$ is not in the positive cone, but
$x(1+x)$ is, and $1+x$ is an order unit. Thus order unit cancellation
fails here. On the other hand, in $R_2$, it is easy to check that the
nontrivial order ideals are exactly $(x^j)$ with $j \geq 2$. The only
prime ideal containing a nontrivial order ideal is thus $(x)$, which
although it is not an order ideal, its square is, and so the criterion of
Theorem 1 applies---$R_2$ does satisfy order unit cancellation. 

The difference in the behaviour of these two examples lies in the
behaviour of the maximal order ideals. In $R_1$, it is $(x(1+x))$ which
(as an ideal) is a product of two distinct prime ideals one of which is
not contained in any order ideal (the split case), while for $R_2$, the
maximal order ideal is $(x^2)$ which is contained in just one prime ideal
of which it is the square (the ramified case). 

\subsubtitle Interesting examples. [H85] Let $H$ be a compact group with
identity $e$, and representation ring denoted $\RR (H)$, and $\Arrow
\pi;H.M_d \C$ a $d$-dimensional representation with character $\chi$. By
reducing to a power and possibly factoring out a normal subgroup, we may
assume in what follows that $\pi$ (respectively $\chi$) is {\it
projectively faithful\/}, that is, $\pi^{-1}(\C\cdot \I) = \brcs{e}$
($\Set{h \in H}{|\chi(h)| = \chi(e)} = \brcs{e}$). Multiplication by
$\chi$ is a positive
$\RR(H)$-module endomorphism of $\RR(H)$, and we can take the direct
limit, forming $\RR(H) [\chi^{-1}] = \lim \Arrow \times
\chi;\RR(H).\RR(H)$, with the direct limit ordering. The convex subgroup
generated by $(\chi_0,1)$ the image of the trivial character at the first
level, or equivalently, the identity in $\RR(H) [\chi^{-1}]$ is a
partially ordered ring (in our strong sense), denoted $R_{\chi}$.

When $H =\T^n$, the $n$-torus, $\chi = P$ is simply a Laurent polynomial
in $n$ variables with nonnegative integer coefficients, and we can permit
nonnegative real coefficients by replacing $\RR(H)$ by $\RR(H)\otimes
\R$. The resulting ordered rings, $R_P$, are the main object of study in
[H87A], and the proof that they satisfy order unit cancellation [H87A, II.5]
is rather intricate, but yields the result referred to earlier about
polynomials in several variables. 

If $H$ is a finite group, then $R_{\chi}$ is simple (as an abelian group)
with unique trace (as follow quickly from the Perron-Frobenius theorem
and the assumption of projective faithfulness), so order unit
cancellation holds in these examples as well.

However, if $H = \T \times \Z_2$, the product of the circle with the
two-element group, there is a projectively faithful character $\chi$ \st
$R_{\chi}$ does not satisfy order unit cancellation [H93, Section 5]. 

If $H$ is a simple connected Lie group (i.e., finite centre), then
generically in $\chi$, $R_{\chi}$ satisfies order unit cancellation, but
it is unknown whether this holds for all projectively faithful $\chi$.
For example, let $P$ be the Laurent polynomial obtained by restricting
$\chi$ to a maximal torus. The Weyl group, $W$, acts on $R_P$, and there
is a natural inclusion 
$R_{\chi} \to R_P$ whose image is contained in $R_P^W$, the fixed point
subring under the action of $W$. If all the extreme points of the convex
hull of the weights of all the irreducibles in $\chi$ do not lie
on the boundary of the Weyl chamber, then the map $R_{\chi} \to R_P^W$ is
an isomorphism as partially ordered rings (where $R_P^W$ inherits the
ordering from $R_P$) [H95, Theorem 1.1], and it follows immediately from $R_P$ satisfying
order unit cancellation that $R_{\chi}$  does as well.

In the remaining situations, we do not know whether order unit
cancellation holds. There are plenty of cases wherein the map $R_{\chi}
\to R_P^W$ is an isomorphism of rings but not of ordered rings (ring
isomorphism occurs but not necessarily order isomorphism, e.g., if $H$ is
simply laced and
$\chi$ is a power of an irreducible; more generally, iff some power of
$\chi$ is saturated in a certain sense  ([H95, Proposition 2.7] using a
result of O'Brien [O'B]). There are even cases wherein the map is not  
onto. But in all cases, the order units are the same, in the sense that
for $u$ in
$R_{\chi}$, $u$ is an order unit for $R_{\chi}$ if and only if its image
in  $R_P$ is an order unit therein. 

Power series versions of these constructions have also been studied. 
Let $P = \sum a_n x^n = \sum (P,x^n) x^n$ be a convergent Maclaurin
series with radius of convergence $1$, and all coefficients nonnegative.
Define $R_P$ to be what we would get from the direct limit construction
(modified for functions analytic on the disk),
$$ 
R_P = \Set{f/P^k}{\text{there exists positive integer $N$ \st for all $n$,
$|(f,x^n)|
\leq N (P^k, x^n)$}}.
$$
 This is a partially ordered ring. If $P$ is continuous at $1$ and the
coefficients behave reasonably smoothly (no large gaps, etc; e.g., $P =
\sum x^n/n^2$), then $R_P$ is noetherian, satisfies order unit
cancellation for very easy reasons (the only relevant prime ideal
consists of the functions vanishing at zero). On
the other hand, if $P = (a+bx)/(1-x^2)$ with $a$ and $b$ positive, and $a
\neq b$, then $R_P$ does not satisfy order unit cancellation [H03, Proposition B.2 and Corollary B.6] (and
there are many other examples of this kind). However, it would be
interest to decide whether $R_P$ satisfies order unit cancellation when
$a=b$, i.e., $P = (1-x)^{-1}$. In all these cases, provided the
coefficients behave reasonably, the pure point evaluations at $[0,1)$ 
are dense in the pure trace space, but when $P(1)$ does not exist, just
which compactification occurs is unknown. 

\subsubtitle Matrix-related examples. Yet another class of interesting
partially ordered rings (not all, but generically, commutative
non-noetherian domains) appears in [H09]. Let $M$ be a square matrix
whose entries are real polynomials in several variables with no negative
entries, form the direct limit (repeating $M$), and let $S$ be the
centralizer of the induced automorphism in the direct limit module, and
finally, let $R$ be the convex subgroup thereof generated by the
identity. This is a partially ordered ring (not always commutative, or a
domain, or noetherian, but having $1$ as an order unit). The example
mentioned above arising from $\T \times \Z_2$ also appears in this
context, so order unit cancellation need not hold---but it should hold
when $R$ has no zero divisors. 

\subtitle Some positivity criteria

Here we give criteria for positivity,
mostly recalled from earlier work. For $s$ in $R$, let $S(s) = \Set{r\in
R}{r s
\geq 0}$; then
$sS(s)$ is the set of products $rs$ \st $rs \geq 0$, and in particular is
a subset of $R^+$. Define $S^+(s) = S(s) \cap R^+$. 

If $A$ is a subset of $R$, then as usual $(A)$ denotes the ideal generated by $R$. If $T$ is a subset of $R^+$, the {\it order ideal generated by $T$} is 
$$
<T>:= \Set{r \in R}{\text{$\exists$ positive integer $M$ and finite $\brcs{t_i}\subseteq T $ \st $-M\sum t_i \leq r \leq M\sum t_i$}}.
$$
This is the smallest order ideal containing $T$, and when $1$ is an order unit, $<T>$ is an order ideal. (In general, there is no satisfactory way of defining $<T>$, the order ideal generated by $T$, if $T$ has nonpositive elements.)

The approach in the
following is mostly based on that of [H87A, Section II].

\Lem Lemma 2. Let $R$ be a commutative ordered ring with $1$ as an order
unit, and let $s$ be an element of $R$. Suppose that with $J:= \oid (sS^+
{(s)})$, the element $s+J$ is in $(R/J)^+$ (this includes the case that
$s+J = 0$, i.e., $s \in J$). Then $s \in R^+$.

\Pf  We may write $s = p + q$ where $ p \in R^+$ and $q \in J$. There
exists a finite set $r_i$ of elements of $S^+ (s)$ \st  both $\pm q \leq
\sum sr_i$ (the $r_i$s may be repeated, hence we don't need an integer
multiple appearing on the right). Let
$J_0 := \oid (\brcs{sr_i})$, and let $J_1 := \oid(\brcs{sr_i}, p)$ (there
are no problems with the order ideals generated by these sets, since all
the members thereof are in $R^+$). Then $s \in J_1$, and since $J_1$ is
generated as an order ideal by a finite set of elements of the positive
cone of $R$, $J_1$ admits an order unit. We will show that $\gamma(s) >
0$ for all pure (normalized) traces of $J_1$, from which it follows that
$s$ is an order unit of $J_1$ and thus is in $R^+$.

By [H85, Lemma 2.1(c)], $\gamma$ is of either of two forms. Either (i) there
exists a pure trace $L$ of $R$ \st $\gamma = \alpha L\left|_R\right.$
(i.e., up to renormalization by the value of $L$ at the order unit of
$J_1$, $\gamma$  is the restriction of $L$ to $J_1$) for some real $\alpha
> 0$, or (ii) there exists a pure trace $L$ of $R$ \st $L(J_1) = 0$ and
for all $j$ in $J_1$ and $r$ in $R$, $\gamma (j r) = \gamma (j) L(r)$.

Suppose $\gamma (s) = 0$. If (i) holds, then $\gamma(sr_i) = \alpha
L(sr_i) = \alpha L(s)L(r_i) = 0$, whereas if (ii)  holds, then
$\gamma(sr_i) = \gamma(s) L(r_i) = 0$. Hence in either case, $\gamma
(sr_i ) = 0$, and thus $\gamma (q)  = 0$, whence $\gamma (p) = 0$. But
this means $\gamma$ vanishes on the order ideal generated by $\brcs{sr_i}
\cup \brcs{p}$, that is, on $J_1$, a contradiction. Hence $\gamma (s) \neq
0$.

If $\gamma (s ) < 0$, then since $sr_i \in R^+$ and $ 0 \leq \gamma
(sr_i) = \gamma (s) L(r_i)$ (in either case (i) or (ii)) and $r_i \geq
0$, we deduce $L(r_i ) = 0$ for all $i$. Therefore $\gamma (sr_i) = 0$,
which again implies $\gamma (q) = 0$, and thus $\gamma (s) = \gamma (p)
\geq 0$ a contradiction.
\qed

This can be extended. But for now, we note a consequence: if $u$ is an
order unit and $us \geq 0$, then $u$ is in $S^+ (s)$, and we may form the
order ideal $I$ generated by $us$. If $s$ belongs to $I$ (or if $s + I
\in (R/I)^+$), then $s \in J$ (defined in the lemma), so that $s \in
R^+$. It is also true if instead of $s\in I$, we have  $-s + I \in
(R/I)^+$, by an almost identical proof (but there is no particular use
apparent for this weird result). 

Now we can give a proof of Theorem 1. 

\subsubtitle Proof of Theorem 1. Suppose $u$ is an order unit of $R$ and
$a$ is an element thereof \st
$ua \in R^+$. Let $I$ be the order ideal generated by $ua$; then
sufficient for $a \in R^+$ is that $a \in I$.
 
By the Lasker-Noether theorem, there exist meet-irreducible (hence
primary) ideals $Q_j$ with associated minimal (over $I$) ideals $P_j$ \st
$I = \cap Q_j$. In particular, for each $j$, $ua \in Q_j$. Moreover, no
power of $u$ can belong to any proper order ideal, so if $u^m$ belonged to
$Q_j$, then $u^{mk} $ would belong to $P_j^k$ (where $k \equiv k(P_j)$ as
hypothesized), whence it would belong to the order ideal $J$, a
contradiction. Since meet-irreducible ideals are primary (zero divisors
are nilpotent), $ua \in Q_j$ and $u^n \notin Q_j$ for all $n$ imply that
$a \in Q_j$, whence $a \in I$. Now the observation above applies. 
\qed

(Thus far at least, all the noetherian domains which fail to satisfy the
hypotheses in Theorem 1 also fail to satisfy order unit cancellation.) 
 
We will show that all the rings of the form $\R[K]$ with $K $ a compact
convex  polyhedral body in Euclidean space, satisfy order unit
cancellation, while some of their close relatives do not. Recall that the positive cone is generated multiplicatively and additively by the linear forms $\beta_i = a_{0i} + \sum a_{ji} x_j$ with real coefficients $a_{ji}$, \st $\beta_i^{-1}0 \cap K = F_i$ run over all the facets (maximal proper
faces) of $K$ with $\beta_i|K \geq 0$.
 
Let $K$ be a compact convex polyhedron with interior in $\R^n$, and let
$R$ be the partially ordered ring $\R[K]$. Let $I$ be an order ideal
therein, and let $Z(I)$ be zero set of $I$; since order ideals are
generated by monomials in linear forms, the zero set (in $\C^n$) is an
intersection of a union of hyperflats, i.e., a union of flats, and
moreover, is just the complexification of the set of real zeros (again a
union of flats, this time in $\R^n$). Now $Z(I) \cap K$ is a union of
faces in $K$, and our immediate aim is to show that it is Zariski dense
in $Z(I)$, that is, for
$f$ in $R$, $f|Z(I)\cap K \equiv 0$ entails $f| Z(I) \equiv 0$. Without loss
of generality, we may assume (since order ideals are generated by
monomials in the $\beta$s), that the zero set notation $Z(I)$ refers to
real zeros (i.e., $Z(I)$ can be reduced to $Z(I ) \cap \R^n$, which we
relabel as $Z(I)$.
 
To see why this is not entirely trivial, consider a trapezoid, $K$, which
is not a parallelogram, and let $F_1$ and $F_3$ be the two non-contiguous
but non-parallel edges. Their intersection is empty, but the intersection
of the their affine hulls (i.e., the real flats they generate) is a single
point, outside $K$. If $\beta_1$ and $\beta_3$ are the corresponding
linear forms (inward normals) exposing the edges, then the ideal $I$
generated by $\brcs{\beta_1, \beta_3}$ is generated by two of the
generators of $\R[K]^+$, is a maximal ideal, its zero set consists of a
singleton (the point where the flats meet) outside $K$, and the
intersection of the zero set with $K$ is empty.
 
Of course, in this case $\beta_1 + \beta_3$ is an order unit of $R $, so
any order ideal containing it would have to be $R$ itself---in particular,
even though $I$ is generated by monomials in the generating linear forms,
$Z(I) \cap K$ is not Zariski dense in $Z(I)$. So the order ideal
hypothesis is crucial.

Another similar example arises from a pyramid with square base. If we take two non-contiguous triangular faces abutting the apex $v$, their intersection is just a singleton consisting of the apex; however, the intersection of their affine hulls is a line, and in this case, the linear form $\beta_1 + \beta_3$ is a (relative) order unit for the maximal ideal (and in this case, order ideal) consisting of functions vanishing at $v$.

\Lem Proposition 3. Let $G$ be a proper face of $K$. Then
$$
I(G) := \Set{f \in R}{f|G \equiv 0}
$$
is an order ideal, and the quotient $R/I(G)$ is order isomorphic to
$\R[G]$ where
$G$ is embedded in $\R^{n-d}$. If
$S =
\brcs{\beta_i}$ is a linearly independent collection  of cardinality
$n-\dim G$ \st
$\beta_i|G
\equiv 0$, then $S$ is a generating set for $I$ as an ideal, and $\sum
\beta_i$ is an order unit of $I$ (as a partially ordered vector space). 

\Pf First, we note that if $G \subset F_i$ where $F_i$ is a facet exposed
by the (positive in $R$) linear form $\beta_i$, then $\beta_i$ belongs to
$I(G)$. First, we show that  the set of all these $\beta_i$ generate
$I(G)$ as an ideal. If $\dim G = 0$, then $G$ is a vertex (extreme
point), so after a change of variables and relabelling, we may assume
that it is the origin and 
$\beta_i = x_i$ for $1 \leq i \leq n$, and there may be subsequent
$\beta_i$ exposing the origin. Now $(x_i)$ is a maximal ideal contained
in $I(G)$, and therefore must equal $I(G)$. 

More generally, let $g = \dim G$; we may assume that the origin is in
the relative interior of $G$, and we find a collection of facets 
$F_1, F_2,
\dots, F_{n-g}$ \st $\cap F_i = G$, and $F_i$ is the zero set of $x_i$
(after a change of variables). Then it is easy to check that $I(G) =
(x_i)$.

In particular, $I(G)$ is generated as an ideal by positive elements,
hence is directed (equal to the set of differences of its positive
elements). It is also convex, that is, if $0 \leq h \leq f$ (meaning \wrt
$R^+$) and $f \in I(G)$ then  $0 \leq h|G \leq f|G \equiv 0$, so $h|G
\equiv 0$, and thus $h \in I(G)$. Thus $I(G)$ is an order ideal. 

With the origin in the relative interior of $G$, let the copy of Euclidean
space for realizing the quotient algebra be the vector space spanned by
$G$. It is easy to check that the ordering on the quotient is exactly
that spanned by the inward normals to relative facets of $G$, which is
the ordering that we put on $R[G]$. 

The linear independence argument is straightforward.

For any order ideal
$J$, if
$\brcs{a_{\alpha}}$ is a finite generating set consisting of elements of
$J^+ = J \cap R^+$, it is trivial to verify that $\sum a_{\alpha}$ is a
(relative) order unit of
$J$ (this requires the fact that $1$ is an order unit of $R$, which is
part of the definition of partially ordered ring that we are using here). 
\qed

\Lem Proposition 4. Let $G$ be a proper face of $K$, and let $\beta$ be a
linear form that is in $R^+$ and for which $\beta^{-1}0 \cap K =
G$.{\par}
\item{(i)}If
$\gamma$ is any other linear form in $R^+$ \st $\gamma| G \equiv 0$, then
there exists a positive integer $M $ \st $M\beta - \gamma \in R^+$.{\par}
\item{(ii)} Suppose $\brcs{F_i}_{i \in S}$ runs over {\it all\/} facets of
$K$ that contain $G$, and let $\beta_i$ denote the corresponding linear
forms. Then there exist real $a_i > 0$ \st $\beta = \sum_S a_i \beta_i$.
{\par}
\item{(iii)} $\beta$ is an  order unit of the order ideal $I(G)$.

\Pf Without loss of generality, we may assume that either $G$ is the
singleton consisting of the origin, or the origin is in the relative
interior of $G$. Translate $G$ by an element  $v$ so that $v+G$ is in the
interior of $K$, and draw a hyperflat $L$ through $v+G$ that misses $G$;
we can make the norm of
$v$ sufficiently small so that $E = L \cap K$ only hits facets that abut
$G$. 

Now since the origin is in $G$, $\beta$ and $\gamma$ are linear, and kill
the vector space spanned by $G$. Since $E$ is compact, and $E$ misses
$G$, we have that $\inf \beta|E = t > 0$, and similarly, $\sup \gamma|E =
u > 0$ exists, so that on $E$, we have $\beta- (t/2u)\gamma|E \geq 0$.
Since $\beta$ and $\gamma$ are linear, so is $\rho: = \beta - (t/2u)
\gamma$, and thus $\rho$ is positive on $\cup_{\lambda > 0} \lambda E$,
and this includes $K \setminus G$. Since $\rho $ vanishes on $G$, we have
that $\rho|K \geq 0$. By [H88, Lemma 1.1], $\rho \in R^+$. We can thus set $M =
2u/t$. This yields (i). 

Form $\gamma = \sum_S \beta_i$, the $\beta_i$ the linear forms
corresponding to all the facets containing $G$, and find $M$ \st $\rho:= M
\beta  - \gamma \in R^+$. By [H88, Lemma 1.1], we can write $\rho = \sum b_i
\beta_i + \sum B_{\alpha} \beta_{\alpha}$ where the $\beta_{\alpha}$ run
over the remaining linear forms (whose facets don't contain $G$), and the
$a$s are all nonnegative. Evaluating at any point $v$ in the relative
interior of
$G$, we deduce
$0 =
\sum B_{\alpha} \beta_{\alpha} (v)$; but all $\beta_{\alpha}$ are
strictly positive on the relative interior of $G$. Hence all the
$B_{\alpha} = 0$. Part (ii) follows immediately (since $b_i \geq 0$). 

Part (iii) is a consequence of part (i) and the previous result.
\qed

\def\rad{\text{rad}\,}

Let $I$ be an  ideal of the commutative (unital) ring $R$. As usual, we
define the radical of $I$, $\rad I := \Set{r \in R}{\text{there exists
$n$ \st $r^n \in I$}}$. As is well known, $\rad I$ is an ideal, and the
primes minimal over $I$ are the same as those minimal over $\rad I$. If
$R$ is now a partially ordered ring (i.e., $1$ is an order unit of
$R$), and $I$ is an order ideal, $\rad I$
need not be an order ideal, but at
least it is convex, that is, if $0 \leq a \leq b \in \rad I$, then $a \in
I$---to see this, there exists $n$ \st $b^n \in I$, so $0 \leq a^n \leq
b^n \in I$, and since $I$ is convex, $a^n \in I$, and thus $a \in \rad
I$. 

This means we can define an order ideal of $R$, denoted $\rad^+ I$, via
$\rad^+ I = (\rad I)^+ - (\rad I)^+$ (where $(\rad I)^+ = (\rad I) \cap
R^+$. It is routine to verify, in view of the convexity result of the
previous paragraph, that
$\rad^+ I$ is an order ideal, and the inclusions $I \subseteq \rad^+ I
\subseteq \rad I$ guarantee that the minimal prime ideals over each of
these ideals are identical. Thus for the purpose of showing $Z(I) \cap K
$ is Zariski dense in $Z(I)$, we may assume that $I = \rad^+ I$, that is,
the quotient ordered ring $R/I$ has no positive nilpotents. 

We define a {\it flat\/} in $\R^n$ to be a translate of a subspace, and a
{\it hyperflat\/} is a flat of dimension $n-1$. If $\beta^w = \prod
\beta_i^{w(i)}$ is a monomial in the $\beta_i$, then the zero set of
$\beta^w$ is a union of hyperflats, each of whose intersection with $K$
is a facet, $F_i = \beta_i^{-1}0 \cap K$. Since order ideals are
generated (as ideals) by monomials in the $\beta_i$, $Z(I)$ will thus be
an intersection of union of hyperflats, in particular, it will be a union
of flats, which we may write as
$Z(I) = \cup E_j$, where the $E_j$ are pairwise incomparable. As we saw
in an earlier example (without the order ideal hypothesis, merely
generation by monomials), it could happen that for some $j$, $E_j \cap K
$ has smaller dimension than $E_j$ or could even be empty. 

However, each $E_j \cap K$ (if nonempty) is a face, call it $G_j$, of
$K$, and we have $K \cap\(\cup E_j\) = \cup G_j$. For a face $G$ of $K$,
denote by $\ti G$ the affine hull of $G$, that is, the flat it generates,
alternatively, 
$\Set{\sum a_s g_s}{a_s \in \R, \ \sum a_s = 1, \text{ and } g_s \in G}$.

Write, for each $j$, $E_j = \cap_t \ti F_{j,t}$, where the $F_{j,t}$ vary
over the facets that arise from the monomials that yield the $E_j$---that
is, there exist monomials $\beta^w$ in $I$ \st $\beta_{j,t}$ appear (with
nonzero exponent) in $\beta^w$. 

Form the face $H_j = E_j \cap K = \cap_j F_{j,t}$, expressed as an
intersection of facets. Let $\FF_j$  denote the set of {\it all\/} facets
of $K$ that contain $H_j$. Form $\beta^{(j)} = \sum_t \beta_{j,t}$; this
is a linear form and a positive element of $R$ for which $\beta^{-1}0
\cap K = H_j$. By the earlier lemma, we can decompose $\beta^{(j)} =
\sum_{\beta_{i,j}^{-1}0 \cap K \in \FF_j} a_{i,j} \beta_{i,j}$ with
$a_{i,j}$ being {\it positive\/} real numbers.

Set $\gamma = \prod_j \beta^{(j)}$; since $\cap \beta_{j,t}^{-1}0 = E_j$,
we have from the real  (here a trivial consequence of the usual) version
of the Nullstellensatz, that some power of $\gamma$ belongs to $I$; as
$\gamma
\in R^+$, our reduction to $I = \rad^+ I$ entails that $\gamma \in I$.

For each $j$, select a $\beta_{i,j}$ (for which $\beta_{i,j}^{-1}0 \cap K
\in \FF_j$), and take the product over  $j$; take all such products.
From $\gamma \in I$ and $I$ being an order ideal, it easily follows that
all such products belong to $I$. The zero set of each such product is a
union of hyperflats, and it is easy to check that their intersection
(over all the products) is $\cap \ti H_j$. Hence $Z(I) = \cap \ti H_j$,
and $Z(I)\cap K= H_j$, as desired.

\Lem Theorem 5. Let $K$ be a compact convex polyhedral body in $\R^n$.
Then the partially ordered ring $R = \R[K]$ satisfies order unit
cancellation. 

\Pf Let $I$ be an order ideal, with $Z(I) \cap K$ the union of faces,
$\cup G_j$. By above, $Z(I) = \cup \ti G_j$, and by the Nullstellensatz,
any minimal prime ideal over $I$ is minimal over the ideal of functions
vanishing on $Z(I)$ (which we do not assume is an order ideal). For any
face
$G$,
$I(G)$ is an order ideal, and from the generating elements, it is
immediate that its zero set is $\ti G$. Moreover $I(G)$ is a prime ideal.
Hence the minimal prime ideals sitting over $I$  are contained in at
least one of the $I(G_j)$, and the criterion of Theorem 1 is now satisfied. 
\qed

\subtitle More exciting examples

Close relatives of the $\R[K]$ were described in [H88, Example I.5]. 
Take $R = \R [x,y]$ with the partial ordering generated additively and
multiplicatively by the set $\brcs{\R^+, x,y, 1 - (x+3/5)^2 - (y+3/5)^2:=
\alpha}$; it was shown [op\. cit.] that $1$
is an order unit, and obviously $R$ is unperforated.

The pure traces correspond to points in the disk $(x+3/5)^2 +(y+3/5)^2
\leq 1$ in the positive quadrant (since $x$ and $y$ must also be sent to
nonnegative reals). Then  $\beta:= 1/5 - y$ is nonnegative as a function
on the pure trace space, vanishing only at the point $(0,1/5)$, and $ u= y
+ 7/5$ doesn't vanish at all on the pure trace space, so is an order unit
of $R$. Now
$$ \eqalign{
\(\frac15 - y\)\(y + \frac 75\) &= \frac 7{25} - \frac 65 y - y^2 \cr & =
1 - \(x+\frac35\)^2 - \(y+\frac 35\)^2 + x^2 + \frac 65 x = \alpha + x^2 
+
\frac 65 x .\cr  
}$$ In particular, $u \beta \in R^+$. However, it is easy
to show (now) that
$\beta$ itself is not in the positive cone. Any element of the positive
cone is a sum of products of powers of $x, y, \alpha$. The zero set (in
$\R^2$) of these elements are respectively the $y$-axis, the $x$-axis, and
the circle $\alpha =0$. Any sum of products of the generators of the positive cone has a zero
set whose intersection with the union can only be one of the following
types: the empty set, the two-point sets $\brcs{(0,1/5), (0,-7/5)}$ and
$\brcs{(1/5,0),(-7/5,0)}$, the singleton $(0,0)$, the zero sets of $x$,
$y$, and $\alpha$, and all possible unions of these. In particular, if
$(0,1/5)$ belongs to the zero set of a positive element, then so must
$(0,-7/5)$. Hence $\beta = 1/5 - y$ is not positive.

The particular property which helps to explain this is the maximal order
ideal $I = (\brcs{x,\alpha})$ (i.e., generated as an ideal). To see that this {\it
is\/} an order ideal, form $r:= \alpha + 6/5 x + x^2$, and observe that it
is in the ideal, is positive, and vanishes at the trace evaluation
$(0,1/5)$; the latter implies it is not an order unit, so the order ideal
it generates, $J = \langle r \rangle$, is proper. Moreover, $x$ and $\alpha$
belong to $J$, and since order ideals are ideals (when $1$ is an order
unit), so $I \subseteq J$.

On the other hand, the zero set of the ideal $I$ has just two real points
(and these are the only two complex zeroes as well), so $I$ has
codimension $2$, and $R/I$ is spanned by $\brcs{1/5-y + I, 7/5 + y + I}$.
Any proper order ideal must have codimension at least two (by the remark
about zeros of positive elements), so that $I = J$ and $J$ is a maximal
order ideal. Similarly, the ideal $(y,\alpha)$ is a maximal order ideal of
codimension two. The only other maximal order ideal is $(x,y)$ which of
course is a maximal ideal.

Now $I$ is the intersection of the two  maximal ideals $(x,1/5 - y)$ and
$(x,7/5+y)$, and this is what makes the counter-example
work---multiplication by $7/5 + y$ brings the element $1/5-y$ into the
order ideal $I$---if $1/5 - y$ were positive, then any order ideal
containing $u(1/5 -y)$ would also contain
$1/5-y$, yielding an alternative proof of non-positivity.

But we have slightly more in this example. It is easy to check that the
principal ideals $(x)$, $(y)$, and $(\alpha)$ are all order ideals (and of
course prime)---this just uses the fact that if a polynomial in two
variables vanishes on a relatively open subset of a one-dimensional
irreducible algebraic curve, then it vanishes on the whole curve. We may
thus form the quotient ordered ring $Q = R/(\alpha)$ with the quotient
ordering. As a ring, $Q = \R [\cos \theta, \sin \theta]$, and is the
integral closure of $\R [\cos \theta]$ (read $ \R[X]$) in the quadratic
extension field $\R(\cos \theta)[\sin \theta]$ (read
$\R(X)[\sqrt{1-X^2}]$), an easy exercise. The partial ordering is obtained
from the substitution $(x,y) \mapsto( \cos \theta - 3/5, \sin \theta -
3/5)$.

Being integrally closed, finitely generated, and of Krull dimension one,
$Q$ is a Dedekind domain; the maximal ideals correspond precisely to the
points of the circle $\alpha = 0$, and rotation induces an isomorphism
among all the maximal ideals, and again it is easy to check that the class
group has order two. As a partially ordered ring, $1$ is an order unit,
and the trace space is just the portion of the circle in the positive
quadrant. Even here, we obtain a counter-example to order unit
cancellation. Let $B = \beta +(\alpha)$ be the image of $1/5-y$ in $Q$; it
vanishes at a surviving trace (evaluate at $(0,1/5)$, so could not be an
order unit; therefore if it were positive, it would generated a proper
order ideal, call it $J$. If $\Arrow \pi; R .Q$ denotes the quotient map,
then $\pi^{-1}J$ would be an order ideal in $R$, and this would contain
the  ideal $(\alpha, 1/5-y)$, which is not contained in any maximal order
ideal (we have listed the only three), a contradiction. The equation
$(u+ I)B \in Q^+$ holds directly, and of course $\pi (u)$ is an order
unit.

So we have a counter-example in the ring $Q$. We can generalize this
somewhat. We have two maximal ideals sitting over the order ideal $I +
(\alpha)$, one of which, $M_1$, is evaluation at $y = 1/5$, a trace, and
the other, $M_2$ is evaluation at $y = -7/5$, and we also have $M_1 M_2
\subseteq I +( \alpha)$. Let $X$ be the pure trace space; then there
exists $r$ in $M_2$ \st $r(t) \neq 0$ from $X$ being a compact subset of
the maximal ideal space (identify the pure traces with their kernels)
missing $M_2$. Then $\pm r$ is necessarily an order unit (since the pure
trace space in this case is connected, the sign of $r(t)$ is
constant---can fix this even if pure trace space is not connected, just
uses density of the polynomial ring in continuous function space: if $R
\subset C(X)$ separating points and $p$ is a point of $X$ and $Y$ is a
compact subset of $X$ missing $p$, then there exists $r$ in $R$ \st $r(p)
= 0$ and $r|Y$ never zero---take sum of squares). Then $r$ is an order
unit, and maps $M_1$ into $I$.

By adjusting centres and radii of the circle, we can make the arc
approximate as closely as we like the line segment joining $(0,1)$
 to $(1,0)$, and in all these cases order unit cancellation fails. In
some sense, these are perturbations of $\R[K] $ where $K$ is the
standard triangle, which of course does satisfy order unit cancellation.

\subtitle Other positivity results

Some of these might eventually be useful.

\Lem Proposition 6. Suppose that $R$ is a partially ordered ring \st every
order ideal which is meet-irreducible as an order ideal is primary as an
ideal (in particular, if it is meet-irreducible as an ideal and the ring
is noetherian). Then $R$ satisfies order unit cancellation.

\Pf Suppose $u$ is an order unit of $R$, $a$ is an element of $R$, and
$ua \in R^+$. Form $I$, the order ideal generated by the positive element
$ua$. If $a \in I$, we are done, by Lemma 2. Otherwise, assume $a$ does
not belong to $I$. An easy Zorn's lemma argument yields a proper  order
ideal, $J$, maximal \wrt containing $I$ but not $a$. Any larger order
ideal would contain $a$, so $J$ is meet-irreducible as an order ideal. By
hypothesis, $J$ is primary as an ideal. Since all powers of $u$ are order
units and order units cannot belong to a proper order ideal, it follows
that $u^n \notin J$ for all $n$. Since $ua \in J$ and $J$  is primary, it
follows that $a \in J$, a contradiction.
\qed

Let $C(R) = \Set{b \in R^+}{\exists  r \in R^+\setminus \brcs{0}, \exists M
\in \N, \exist \text{ positive invertible } u \text{ \st} r \geq b^M u}$.
For example, if $R^+$ is finitely generated (meaning there is a finite
subset of $R^+$ \st every element of $R^+$ can be expressed as positive
real linear combination of products of elements of the set), then $C(R)$
is nonempty (take the product of the generators). We can weaken the
definition of finite generation to replace the positive real combination
by positive combination with coefficients order units in place of reals.
This includes all the ordered rings of the form $R_P$ studied in [H87a],
e.g., if $P = 1 + x + y$, set $b = xy/P^3$. The rings obtained from power
series in [H96], [H03] (also denoted $R_P$) do not have  finitely generated
positive cones, but $C(R_P)$ are not empty here as well---if $P(0) \neq
0$, take $b = P^{-1}$.

Finally, for $s$ in $R$, let $S(s) = \Set{r\in R}{r s \geq 0}$; then
$sS(s)$ is the set of products $rs$ \st $rs \geq 0$, and in particular is
a subset of $R^+$.

The following is a partial improvement on [H03, Lemma 4.2(b)], which in turn was based on [H95, Corollary 1.3]. For $R$ a partially ordered ring (with $1$ as order unit), the pure trace space (consisting of the positive ring homomorphisms $R \to \R$ and equipped with the point-open topology) is denoted $X(R)$, and is a compact Hausdorff space. 

\Lem Proposition 7. Suppose $X(R)$ is connected and $s$
is an element of $R$ with the following properties{\par}
\item{(a)} $\Set{x \in X}{ x(s) = 0}$ is nowhere dense in $X${\par}
\item{(b)} $\ideal(\id{sS(s)}, S(s)) = R$.{\par}
 \noindent Then one of $
\pm s$ belongs to $R^+$.

\Lem Corollary 8. Suppose $X(R)$ is connected, and $s$ is a
nonzero-divisor of  $R$ whose zero set (as a function on $X$) is nowhere
dense. Suppose in addition, $\Z[1/p] \subset R\cdot 1$ for some positive
integer $ p >1$.  If $\ideal( \id {C(R),sS(s)}, S(s) ) = R$, then one
of  $ \pm s$ belongs to $R^+$.

\noindent {\it Proof of Corollary 8\/} (from Proposition 7). Suppose $rs
\geq 0$; as $rs \neq 0$, for any $b$ in $C(R)$, there exists a positive
integer $M$ and a positive invertible $u$ in $R$ \st $rs \geq b^M u$.

There exists an integer $n$ \st $s \geq -n$. Now consider for any
multiple of the identity $\epsilon\cdot 1$,
$$
(r + \epsilon b^M u) s = b^M u + \epsilon b^M u s \geq b^M u ( 1 - n
\epsilon).
$$
 Select $\epsilon  < 1/n$ (all that is required  is that the
image of the invertible constants in $R$ be dense in the reals). Thus $r
+ \epsilon b^M u$ belongs to $S(s)$, so that $b^M$ belongs to
$\ideal(S(s))$. Hence for each $b$ in $C(R)$, there exists $M$ (depending
on $b$) \st $b^M \in \ideal(\oid ( s S(s)),S(s))$.

Since $\ideal( \oid (C(R)),S(s) ) = R$, there exist $B_1, \dots , B_N$
in $\oid({C(R)})$ and $r_i$ in $S(s)$ together with $t_i$ in $R$ \st $ 1 =
\sum B_i + \sum r_i t_i$ (an order ideal in an ordered ring with $1$ as
order unit is automatically an ideal). Write each $B_i = d_{i1} - d_{i2}$
with $d_{ij} \geq 0$, so there exists $K$ \st $d_{ij} \leq K \sum b_l$
for some finite collection of $b_l$ in $C(R)$. Hence there exist (large)
$L$ \st $d_{ij}^L \in  \ideal(\oid ( s S(s)),S(s))$, and even larger $J$
\st $B_i^{J}$ also belongs, and it easily follows that $ \ideal(\oid (s
S(s)),S(s))$ is improper, so the preceding result applies.
\qed

\noindent {\it Proof of Proposition 8.}  By [H03, 4.2(B)], it suffices to
show $S(s)$ generates the improper ideal. As in the argument in the proof
of the corollary, if $r
\in S(s)$, then $(r + \epsilon rs) s \geq rs (1 - M\epsilon)$, so if
$\epsilon $ is chosen sufficiently small, $r + \epsilon rs$ belongs to
$S(s)$, and thus $rs \in \ideal (S(s))$. If $0 \leq c_{t1},c_{t2} \leq K
\sum r_i s$ for a finite selection of $r_i$ in $S(s)$ and $c_{tj}$ in 
the positive cone of $\id{sS(s)}$ and $1 = \sum (c_{t1} -c_{t2}) + \sum
r_i s w_i$, then each $r_is $ belongs to $\ideal(S(s))$. Moreover, since
$S(s)$ is closed under sums, the term $\sum r_i s$ belongs to
$\ideal(S(s))$, and now the preceding argument works for each $c_{tj}$
(with appropriately small $\epsilon$). Hence the $c_{tj}$ belong to
$\ideal (S(s))$, so the latter is improper.
\qed

Unfortunately, this is not enough to get the order unit cancellation
results of [H87A, section 2]. With $P = 1 + x +y$ and $b = xy/P^3$ (or anything similar;
it must vanish on the boundary of the Newton polytope). For example, $Y
+3 = y/P + 3$ is an order unit and $\ideal(b, Y+3)$ is proper (both terms
vanish when $Y = -3$ and $x = 0$) (in this case, $\id{C(R)} = bR$).

\subtitle Sums of order ideals

Consider the condition on a partially ordered ring $R$ that (finite) sums
of order ideals are order ideals
This is a consequence of (Reisz) interpolation. Rather remarkably, for
the ordered rings $\R[K]$, it actually implies interpolation, and
therefore that up to affine equivalence, $K$ is  a product of simplices.

The intersection between the two class of ordered rings, the
collection of $\R[K]$ on the one hand, and of $R_P$ (with $P$ a polynomial
in several variables with only nonnegative coefficients) is very tiny. The
ordered rings in the latter class satisfy interpolation (in fact, they
are dimension groups, and their direct limit structure follows from their
definition). It was shown in [H88] that if $\R[K] $ satisfies
interpolation, then $K$ is affinely equivalent to a product of simplices,
and in particular, $\R[K]$ is of the form $R_P$ for very special choices
of $P$ (and as a result, positivity of a given element can be determined
explicitly by the main result of [H86]).

An examination of the argument shows that the hypotheses can be reduced
considerably. In fact, it requires hardly any additional effort to obtain
the following.

\Lem Theorem 9. Let $K$ be a compact convex polyhedron with interior in
$\R^n$. The following are equivalent.{\par}
\item{(a)} $\R[K]$ satisfies the Riesz interpolation property{\par}
\item{(b)} sums of order ideals in $\R[K]$ are order ideals {\par}
\item{(c)} any ideal generated by a subset of $\brcs{\beta_i}$ (linear
forms associated to the facets) is an order ideal {\par}
\item{(d)} $K$ is $\text{AGL}(n,\R)$-equivalent to a product of simplices
{\par}
\item{(e)} there exist $n(i)$-dimensional simplices $K_i$  with $\sum
\text{dim} K_i = n$ \st $\R[K]$ is order isomorphic to $\otimes_{\R}
\R[K_i]$ (tensor product ordering).{\par}

\noindent {\it Remark\/} In part (e), if $K_i$ is $n(i)$ dimensional, then
$\R[K_i]$ is order isomorphic  (i.e., as ordered rings) to $R_{P_i}$ with
$P_i = 1 +
\sum_{j=1}^{n(i)} X_{ji}$, via the map  $\beta_j \mapsto  X_{ji}/P$ and
$\beta_0$ (corresponding to the facet that misses the origin) is sent to
$1/P_i$ (see the discussion at the end of [H88]), and is $R = \R[\prod
K_i]$, then $R \iso R_P$ where $P = \prod P_i$.

\Pf By [H88, Lemma II.1(b)], any principal ideal$(\beta^w)$ (a product of
the irredundant linear forms $\beta_i$) is an order ideal. Thus (b)
implies (c), and now only (c) implies (d) requires proof. 

We adapt the arguments of [H88, II.2 and II.3], in order to verify the
criteria of II.5 therein. First, if $\brcs{F_{\alpha}}$ is a family of
facets with $\cap F_{\alpha} = \emptyset$, then we show the intersection
of the corresponding affine spans, 
$\tilde F_{\alpha}$,
 is also empty.
Let $\beta_{\alpha}$ be the corresponding (irredundant) linear forms
exposing $F_{\alpha}$ (and for which $\beta_{\alpha}|K \geq 0$). Then
$\beta:= \sum \beta_{\alpha}$ vanishes nowhere on $K$ so is an order
unit. Therefore the order ideal generated by $\brcs{\beta_{\alpha}}$ is
improper. Since sums of order ideals are order ideals and each
$\beta_{\alpha} R$
is itself an order ideal (a consequence of [H88,II.1(b)]), it follows that
$\brcs{\beta_{\alpha}}$ generates the improper ideal (as an ideal, not
just as an order ideal). In particular, $\beta_{ \alpha}$ have no common
zero in all of $\R^n$ (where $K \subset \R^n$). This implies $\cap
\tilde F_{\alpha}$ is empty.

If in the preceding, we replace facets by faces, we observe that every
face is an intersection of facets, so we obtain, if $\brcs{F_{\alpha}}$ is
a collection of faces with empty intersection, then $\cap
\tilde F_{\alpha}$ is also empty. This verifies condition (ai) of [H88,
II.5].


Now suppose that there is a vertex lying in more than $n$ facets, $F_1$, $F_2$, \dots, $F_{n+1}$ abutting the vertex $v$; let $\beta_i$ be the
corresponding linear forms. Without loss of generality, we may assume
that $v$ is origin, so that $\beta_i$ have zero constant terms, and are
thus linear. Since the intersection of any $n$ of the facets must consist
only of the vertex, we find that any $n$-element subset of
$\brcs{\beta_i}_{i=1}^{n+1}$ is linearly independent, and therefore spans
$\R^n$. Hence we may find reals
$a_i$ \st $\beta_{n+1} = \sum_{1 \leq i \leq n} a_i \beta_i$. All of the
$a_i$ must be nonzero (else linear independence of any subset of size $n$
or less is violated). If all the $a_i$ are positive, then
$\beta_{n+1}^{-1}0 \cap K$ can only be a singleton, rather than a facet,
a contradiction; if all the $a_i$ are negative, $\beta_{n+1}$ is negative
on $K$, again a contradiction.

Hence $\brcs{1,2,\dots,n} $ decomposes as $S \cup T$, where $S =
\Set{i}{a_i > 0}$ and $T = \Set{i}{a_i < 0}$, and both $S$ and $T$ are
nonempty. We thus have
$\sum_T |a_i| \beta_i \leq \sum_S a_i \beta_i$. From [H88] again, each
$\beta_i R$ is an order ideal, and the hypothesis asserts that $\sum_S
\beta_i R$ is an order ideal. Thus (as $|a_i |> 0$ for $i$ in $T$),
$\brcs{\beta_i}_T \subset \sum_S \beta_i R$; in particular, the ideal of
$R$ generated by $\brcs{\beta_i}_{1\leq i \leq n}$ is generated by
$\brcs{\beta_i}_S$, which has fewer than $n$ elements. This is
impossible---$M = (\beta_i)$ is a maximal ideal, and since $M/M^2$ is
$n$-dimensional, $M$ cannot be generated by fewer than $n$ elements.
\qed

In a dimension group, finite intersections of order ideals are order
ideals. I do not know whether this holds for all $\R[K]$, or whether it
imposes  constraints on $K$. The rings of the form $R_P$ also have the
property that (finite) products of order ideals are order ideals (not all
ordered rings that are also dimension groups satisfy this property).
Again, it would be worth investigating what $\R[K]$ having this  property
says about $K$.

It also follows from the argument of (c) implies (d), that $U^{-1}\R[K]$
satisfies interpolation (or its weaker forms) implies $K$ is simplicial
(exactly $n$ facets hitting each vertex). Whether the converse holds is
unknown.


\subtitle References

\long\def\Reff[#1] #2, #3, #4\par{\vskip 1pt \item{[#1]} #2, {\it #3,}
#4\par} {\parindent = 2.5em 

\Reff [EHS] Edward G Effros{, David E Handelman,} and Chao-Liang Shen,
Dimension Groups and Their Affine Representations, 
American Journal of Mathematics 102 (1980) 385--407.

\Reff [GH] K Goodearl and D
Handelman, Rank functions and K${}_0$ of regular rings, J  Pure Appl. Alg 7 (1976) 195--216.

\Reff  [H81] David Handelman, Rings with involution as partially ordered
abelian groups, Rocky Mountain J Math, Volume 11  (1981) 337--382.
 
\Reff [H85] David Handelman, Positive polynomials and product 
type actions of compact groups, Mem Amer Math Soc, 320 (1985) 79 p+ xi.

\Reff [H86] David Handelman, Deciding eventual positivity of polynomials,
Ergodic Theory and Dynamical Systems (1986) 57--79.

\Reff [H87] David Handelman, Extending traces on fixed point C* 
algebras under xerox product type actions of compact Lie groups, Journal
of functional analysis 72 (1987) 44--57.

\Reff [H87A] David Handelman, Positive polynomials{, convex integral 
polytopes and a random walk problem},   (1987), Springer-Verlag  1282,
142p.
 
\Reff [H88] David Handelman,
Representing polynomials by positive linear functions on compact convex
polyhedra, Pacific J Math Volume 132 (1988) 35--62.

\Reff [H90] David Handelman, Effectiveness of an affine invariant for
indecomposable integral polytopes,  Int J Math 4  (1993)
59--88.

\Reff [H93] David Handelman, Representation rings as invariants for compact groups and ratio limit theorems for them,  J Pure Appl Algebra  66  (1990)
165--184.

\Reff [H95] David Handelman, Iterated multiplication of characters of
compact connected Lie groups,  J Algebra  173  (1995) 67--96.
 
\Reff [H96] David Handelman, Eventual positivity for analytic functions, 
Math Ann, 304  (1996) 315--338.
 
\Reff [H03] David Handelman, More eventual positivity for analytic
functions, Canad J Math 55 (2003) 1019--1079
 
\Reff [H09] David Handelman,  Matrices of positive polynomials, Electron.
J Linear Algebra 19 (2009) 2--89.  

\Reff [O'B] Ellen O'Brien,  A saturation property for powers of
irreducible characters, Crelle's J 456 (1994) 151--171. 

}  
 
\vskip 6pt \noindent Mathematics Department, University of Ottawa, Ottawa
ON K1N 6N5, Canada; dehsg\@uottawa.ca

\end

%% file: generic_macros




\font\rm=cmr10 \rm

\font\bf=cmb10
\font\Rm=cmr9 at 11pt
\rm
\font\it=cmsl9 at 10pt
 at 7pt

\font\Rrm=cmr17 at 16pt
   \font\Rm=cmr12 at 11.5pt

\long\def\Pf{\par\noindent {\it Proof.} }
\def\({\left(}
\def\){\right)}
\def\st{such that }
\def\qed{\hfill$\bullet$\vskip 4pt}
\def\quotes#1{{\lq\lq #1\rq\rq}}
\def\brcs#1{\left\{ #1\right\}}

\def\iso{\cong}
\def\wrt{with respect to }
\def\:{\,:}

\def\ti#1{\tilde{#1}}

\def\C{\text{\bf C}}
\def\T{\text{\bf T}}

\def\I{\text{I\,}}

\def\R{\text{\bf R}}
\def\N{\text{\bf N}}
\def\Z{\text{\bf Z}}

\def\Arrow #1;#2.{#1\:#2 \to }

\def\Set#1#2{\brcs{#1 \left|\vphantom{#1 #2} \right.#2}}



\def\Rrr#1,#2{{\Cal J}_{#1,#2}}
\def\slfrac#1#2{{\raise -.07 ex\hbox{$^{#1}$}}\!/\raise .35 ex \hbox{${}_{#2}$}}
\def\ssf #1/#2{\slfrac {#1}{#2}}

\def\pd #1,#2.{\frac {\partial #1}{\partial #2}}

   \long\def\Lem
#1.#2\par{\vskip4pt{\baselineskip=13pt\font\it=cmsl12 at
11.5pt\Rm
   \noindent {\rm \uppercase{#1}} #2\vskip3pt

   }}

\long\def\Title #1\par {\noindent{\Rrm #1}\vskip 9pt}

 \long\def\SubT #1.{\noindent {\it #1\/} } 
 
 \long\def\SecT
#1\par{\vskip 3pt \noindent {\bf #1}\vglue1pt
   \noindent}

\long\def\subtitle #1.{\vskip 2pt \noindent {\it #1}}

\long\def\Rmk#1\par{\vskip 1pt \noindent {\it
Remark.} #1\vskip2pt}


%% file: papermacros_entirefunctions
\scrollmode\NoBlackBoxes
\magnification=1100
\long\def\Abstract #1\par%
{\vskip .2 true cm{\leftskip 1 true in \rightskip 1 true in \font\rm=cmr8 \rm
\baselineskip=1pt \font\it=cmsl8 \font\bf=cmb8
\parindent=0em {\bf Abstract} #1

}}
\comment
\font\rm=Times at 10pt

\font\bf=TimesB
\font\Rm=Times at 11pt
\rm
\font\it=TimesI at 10pt
\endcomment

\long\def\Pf{\par\noindent {\it Proof.} }
\def\({\left(}
\def\){\right)}
\def\st{such that }
\def\qed{\hfill$\bullet$\vskip 4pt}
\def\quotes#1{{\lq\lq #1\rq\rq}}
\def\brcs#1{\left\{ #1\right\}}
\def\Set#1#2{\brcs{#1 \left|\vphantom{#1 #2} \right.#2}}

\def\C{\text{\bf C}}
\def\T{\text{\bf T}}

\def\I{\text{I\,}}

\def\iso{\cong}
\def\wrt{with respect to }
\def\:{\,:}
\def\Arrow #1;#2.{#1\:#2 \to }


\def\R{\text{\bf R}}
\def\N{\text{\bf N}}
\def\Z{\text{\bf Z}}

\def\Rrr#1,#2{{\Cal J}_{#1,#2}}

\def\slfrac#1#2{{\raise -.07 ex\hbox{$^{#1}$}}\!/\raise .35 ex \hbox{${}_{#2}$}}
\def\ssf #1/#2{\slfrac {#1}{#2}}

\def\pd #1,#2.{\frac {\partial #1}{\partial #2}}


   \long\def\Title #1\par {\noindent{\Rrm #1}\vskip 9pt}
 \long\def\SubT #1.{\noindent {\it #1\/} }   \long\def\SecT
#1\par{\vskip 3pt \noindent {\bf #1}\vglue1pt
   \noindent}
\long\def\subtitle #1.{\vskip 2pt \noindent {\it #1}}

\long\def\Rmk#1\par{\vskip 1pt \noindent {\it
Remark.} #1\vskip2pt}

